\documentclass[12pt,vatola]{article}

\usepackage{graphicx}

\def\C{\centerline}

\def\re#1{\par\hangindent\parindent\indent\llap{#1\enspace}\ignorespaces}

\def\no{\noindent}

\begin{document}

\vskip 15mm

\C{\large\bf On Algebraic Multi-Vector Spaces  }  \vskip 5mm

\C{Linfan Mao} \vskip 3mm \C{\scriptsize (Academy of Mathematics
and System Sciences, Chinese Academy of Sciences, Beijing 100080)}

\vskip 8mm
\begin{minipage}{130mm}
\no{\bf Abstract}: {\small A Smarandache multi-space is a union of
$n$ spaces $A_1,A_2,\cdots ,A_n$ with some additional conditions
holding. Combining Smarandache multi-spaces with linear vector
spaces in classical linear algebra, the conception of multi-vector
spaces is introduced. Some characteristics of a multi-vector space
are obtained in this paper.}

\vskip 2mm \no{\bf Key words:} {\small  vector, multi-space,
multi-vector space, ideal subspace chain.}

 \vskip 2mm \no{{\bf
Classification:} AMS(2000) 15A04, 15A09,15A18}
\end{minipage}

\vskip 8mm

{\bf $1$. Introduction}

\vskip 6mm

The notion of multi-spaces is introduced by Smarandache in $[6]$
under his idea of hybrid mathematics: {\it  combining different
fields into a unifying field}($[7]$), which is defined as follows.

\vskip 4mm

\no{\bf Definition $1.1$} \ {\it For any integer $i, 1\leq i\leq
n$ let $A_i$ be a set with ensemble of law $L_i$, and the
intersection of $k$ sets $A_{i_1},A_{i_2},\cdots , A_{i_k}$ of
them constrains the law $I(A_{i_1},A_{i_2},\cdots , A_{i_k})$.
Then the union of $A_i$, $1\leq i\leq n$

$$\widetilde{A} \ = \ \bigcup\limits_{i=1}^n A_i$$

\no is called a multi-space.}

\vskip 3mm

As we known, a {\it vector space} or {\it linear space} consists
of the following:

($i$) a field $F$ of scalars;

($ii$) a set $V$ of objects, called vectors;

($iii$) an operation, called vector addition, which associates
with each pair of vectors ${\bf a,b}$ in $V$ a vector ${\bf a+ b}$
in $V$, called the sum of ${\bf a}$ and ${\bf b}$, in such a way
that

($1$) addition is commutative, ${\bf a+ b = b+ a}$;

($2$) addition is associative, ${\bf (a+ b)+ c= a+( b+ c)}$;

($3$) there is a unique vector ${\bf 0}$ in $V$, called the zero
vector, such that ${\bf a+ 0= a}$ for all ${\bf a}$ in $V$;

($4$) for each vector ${\bf a}$ in $V$there is a unique vector
${\bf -a}$ in $V$ such that ${\bf a+(-a)= 0}$;

($iv$) an operation ¡°$\cdot$¡±, called scalar multiplication,
which associates with each scalar $k$ in $F$ and a vector ${\bf
a}$ in $V$ a vector $k\cdot {\bf a}$ in $V$, called the product of
$k$ with ${\bf a}$, in such a way that

($1$) $1\cdot {\bf a}={\bf a}$ for every ${\bf a}$ in $V$;

($2$) $(k_1k_2)\cdot {\bf a}=k_1(k_2\cdot {\bf a});$

($3$) $k\cdot ({\bf a}+{\bf b})=k\cdot {\bf a}+k\cdot {\bf b};$

($4$) $(k_1+k_2)\cdot {\bf a}=k_1\cdot {\bf a}+k_2\cdot {\bf a}.$

\no We say that $V$ is a {\it vector space over the field $F$,
denoted by $(V \ ;+,\cdot )$.}

By combining Smarandache multi-spaces with linear spaces, a new
kind of algebraic structure called multi-vector space is found,
which is defined in the following.

\vskip 4mm

\no{\bf Definition $1.2$} \ {\it Let
$\widetilde{V}=\bigcup\limits_{i=1}^k V_i$ be a complete
multi-space with binary operation set $O(\widetilde{V})=\{
(\dot{+}_i,\cdot_i) \ | \ 1\leq i\leq m\}$ and
$\widetilde{F}=\bigcup\limits_{i=1}^k F_i$ a multi-filed space
with double binary operation set $O(\widetilde{F})=\{(+_i,\times_i
)\ | \ 1\leq i\leq k\}$. If for any integers $i,j, \ 1\leq i,
j\leq k$ and $\forall {\bf a,b,c}\in\widetilde{V}$,
$k_1,k_2\in\widetilde{F}$,

$(i)$ $(V_i;\dot{+}_i,\cdot_i)$ is a vector space on $F_i$ with
vector additive $\dot{+}_i$ and scalar multiplication $\cdot_i$;

$(ii)$ $({\bf a}\dot{+}_i{\bf b})\dot{+}_j{\bf c}= {\bf
a}\dot{+}_i({\bf b}\dot{+}_j{\bf c})$;

$(iii)$ $(k_1+_i k_2)\cdot_j{\bf a}=k_1+_i(k_2\cdot_j{\bf a});$

\no if all those operation results exist, then $\widetilde{V}$ is
called a multi-vector space on the multi-filed space
$\widetilde{F}$ with a binary operation set $O(\widetilde{V})$,
denoted by $(\widetilde{V}; \widetilde{F})$.}

\vskip 3mm

For subsets $\widetilde{V}_1\subset\widetilde{V}$ and
$\widetilde{F}_1\subset\widetilde{F}$, if $(\widetilde{V}_1;
\widetilde{F}_1)$ is also a multi-vector space, then call
$(\widetilde{V}_1; \widetilde{F}_1)$ a multi-vector subspace of
$(\widetilde{V}; \widetilde{F})$.

The subject of this paper is to find some characteristics of a
multi-vector space. For terminology and notation not defined here
can be seen in $[1],[3]$ for linear algebraic terminologies and in
$[2],[4]-[11]$ for multi-spaces and logics.

\vskip 8mm

{\bf $2.$ Characteristics of a multi-vector space}

\vskip 5mm

First, we have the following result for multi-vector subspace of a
multi-vector space.

\no{\bf Theorem $2.1$} \ {\it For a multi-vector space
$(\widetilde{V}; \widetilde{F})$,
$\widetilde{V}_1\subset\widetilde{V}$ and
$\widetilde{F}_1\subset\widetilde{F}$, $(\widetilde{V}_1;
\widetilde{F}_1)$ is a multi-vector subspace of $(\widetilde{V};
\widetilde{F})$ if and only if for any vector additive
¡°$\dot{+}$¡±, scalar multiplication ¡°$\cdot$¡± in
$(\widetilde{V}_1; \widetilde{F}_1)$ and $\forall {\bf
a,b}\in\widetilde{V}$, $\forall \alpha\in\widetilde{F}$,}

$$\alpha\cdot {\bf a}\dot{+}{\bf b}\in\widetilde{V}_1$$

\no{\it if their operation result exist.}

\vskip 3mm

{\it Proof} \ Denote by $\widetilde{V}=\bigcup\limits_{i=1}^k V_i,
\widetilde{F}=\bigcup\limits_{i=1}^k F_i$. Notice that
$\widetilde{V}_1=\bigcup\limits_{i=1}^k(\widetilde{V}_1\bigcap
V_i)$. By definition, we know that  $(\widetilde{V}_1;
\widetilde{F}_1)$ is a multi-vector subspace of $(\widetilde{V};
\widetilde{F})$ if and only if for any integer $i, 1\leq i\leq k$,
$(\widetilde{V}_1\bigcap V_i; \dot{+}_i,\cdot_i)$ is a vector
subspace of $(V_i, \dot{+}_i,\cdot_i)$ and $\widetilde{F}_1$ is a
multi-filed subspace of $\widetilde{F}$ or $\widetilde{V}_1\bigcap
V_i=\emptyset$.

According to the criterion for linear subspaces of a linear space
([$3$]), we know that for any integer $i, 1\leq i\leq k$,
$(\widetilde{V}_1\bigcap V_i; \dot{+}_i,\cdot_i)$ is a vector
subspace of $(V_i, \dot{+}_i,\cdot_i)$ if and only if for $\forall
{\bf a, b}\in\widetilde{V}_1\bigcap V_i$, $\alpha\in F_i$,

$$\alpha\cdot_i{\bf a}\dot{+}_i{\bf b}\in\widetilde{V}_1\bigcap V_i.$$

\no That is, for any vector additive $\dot{+}$, scalar
multiplication $\cdot$ in $(\widetilde{V}_1; \widetilde{F}_1)$ and
$\forall {\bf a,b}\in\widetilde{V}$, $\forall
\alpha\in\widetilde{F}$, if $\alpha\cdot {\bf a}\dot{+}{\bf b}$
exists, then $\alpha\cdot {\bf a}\dot{+}{\bf
b}\in\widetilde{V}_1$.\quad\quad $\natural$

\vskip 4mm

\no{\bf Corollary $2.1$} \ {\it Let
$(\widetilde{U};\widetilde{F}_1),(\widetilde{W};\widetilde{F}_2)$
be two multi-vector subspaces of a multi-vector space
$(\widetilde{V};\widetilde{F})$. Then
$(\widetilde{U}\bigcap\widetilde{W};\widetilde{F}_1\bigcap\widetilde{F}_2)$
is a multi-vector space.}

\vskip 3mm

For a multi-vector space $(\widetilde{V};\widetilde{F})$, vectors
${\bf a}_1,{\bf a}_2,\cdots ,{\bf a}_n\in\widetilde{V}$, if there
are scalars $\alpha_1, \alpha_2, \cdots ,\alpha_n\in\widetilde{F}$
such that

$$\alpha_1\cdot_1{\bf a}_1\dot{+}_1\alpha_2\cdot_2{\bf a}_2\dot{+}_2\cdots\dot{+}_{n-1}\alpha_n\cdot_n
{\bf a}_n={\bf 0},$$

\no where ${\bf 0}\in\widetilde{V}$ is an unit under an operation
¡°$+$¡± in $\widetilde{V}$ and $\dot{+}_i, \cdot_i\in
O(\widetilde{V})$, then the vectors ${\bf a}_1,{\bf a}_2,\cdots
,{\bf a}_n$ are said to be {\it linearly dependent}. Otherwise,
${\bf a}_1,{\bf a}_2,\cdots ,{\bf a}_n$ to be {\it linearly
independent}.

Notice that in a multi-vector space, there are two cases for
linearly independent vectors ${\bf a}_1,{\bf a}_2,\cdots ,{\bf
a}_n$:

\vskip 2mm

($i$) \ for any scalars $\alpha_1, \alpha_2, \cdots
,\alpha_n\in\widetilde{F}$, if

$$\alpha_1\cdot_1{\bf a}_1\dot{+}_1\alpha_2\cdot_2{\bf a}_2
\dot{+}_2\cdots\dot{+}_{n-1}\alpha_n\cdot_n {\bf a}_n={\bf 0},$$

\no where ${\bf 0}$ is a unit of $\widetilde{V}$ under an
operation ¡°$+$¡± in $O(\widetilde{V})$, then $\alpha_1=0_{+_1},
\alpha_2=0_{+_2}, \cdots , \alpha_n=0_{+_n}$, where $0_{+_i},
1\leq i\leq n$ are the units under the operation $+_i$ in
$\widetilde{F}$.

($ii$) \ the operation result of $\alpha_1\cdot_1{\bf
a}_1\dot{+}_1\alpha_2\cdot_2{\bf
a}_2\dot{+}_2\cdots\dot{+}_{n-1}\alpha_n\cdot_n {\bf a}_n$ does
not exist.\vskip 2mm

Now for a subset $\widehat{S}\subset\widetilde{V}$, define its
{\it linearly spanning set} $\left<\widehat{S}\right>$ to be

$$\left<\widehat{S}\right> = \{ \ {\bf a} \ | \ {\bf a}=\alpha_1\cdot_1{\bf a}_1\dot{+}_1
\alpha_2\cdot_2{\bf a}_2\dot{+}_2\cdots \in\widetilde{V}, {\bf
a}_i\in\widehat{S}, \alpha_i\in\widetilde{F}, i\geq 1 \}.$$

\no For a multi-vector space $(\widetilde{V};\widetilde{F})$, if
there exists a subset $\widehat{S},
\widehat{S}\subset\widetilde{V}$ such that $\widetilde{V} =
\left<\widehat{S}\right>$, then we say $\widehat{S}$ is a {\it
linearly spanning set} of the multi-vector space $\widetilde{V}$.
If the vectors in a linearly spanning set $\widehat{S}$ of the
multi-vector space $\widetilde{V}$ are linearly independent, then
$\widehat{S}$ is said to be a {\it basis} of $\widetilde{V}$.

\vskip 4mm

\no{\bf Theorem $2.2$} \ {\it Any multi-vector space
$(\widetilde{V};\widetilde{F})$ has a basis.}

\vskip 3mm

{\it Proof} \ Assume $\widetilde{V}=\bigcup\limits_{i=1}^k V_i,
\widetilde{F}=\bigcup\limits_{i=1}^k F_i$ and the basis of the
vector space $(V_i;\dot{+}_i,\cdot_i)$ is $\Delta_i=\{{\bf
a}_{i1},{\bf a}_{i2},\cdots ,{\bf a}_{in_i}\}$, $1\leq i\leq k$.
Define

$$\widehat{\Delta} \ = \ \bigcup\limits_{i=1}^k\Delta_i.$$

\no Then $\widehat{\Delta}$ is a linearly spanning set for
$\widetilde{V}$ by definition.

If vectors in $\widehat{\Delta}$ are linearly independent, then
$\widehat{\Delta}$ is a basis of $\widetilde{V}$. Otherwise,
choose a vector ${\bf b}_1\in\widehat{\Delta}$ and define
$\widehat{\Delta}_1=\widehat{\Delta}\setminus\{{\bf b}_1\}$.

If we have obtained the set $\widehat{\Delta}_s, s\geq 1$ and it
is not a basis, choose a vector ${\bf
b}_{s+1}\in\widehat{\Delta}_s$ and define
$\widehat{\Delta}_{s+1}=\widehat{\Delta}_s\setminus\{{\bf
b}_{s+1}\}$.

If the vectors in $\widehat{\Delta}_{s+1}$ are linearly
independent, then $\widehat{\Delta}_{s+1}$ is a basis of
$\widetilde{V}$. Otherwise, we can define the set
$\widehat{\Delta}_{s+2}$. Continue this process. Notice that for
any integer $i, 1\leq i\leq k$, the vectors in $\Delta_i$ are
linearly independent. Therefore, we can finally get a basis of
$\widetilde{V}$. \quad\quad $\natural$

Now we consider the finite-dimensional multi-vector space. A
multi-vector space $\widetilde{V}$ is {\it finite-dimensional} if
it has a finite basis. By Theorem $2.2$, if for any integer $i,
1\leq i\leq k$, the vector space $(V_i; +_i,\cdot_i)$ is
finite-dimensional, then $(\widetilde{V};\widetilde{F})$ is
finite-dimensional. On the other hand, if there is an integer
$i_0, 1\leq i_0\leq k$, such that the vector space $(V_{i_0};
+_{i_0},\cdot_{i_0})$ is infinite-dimensional, then
$(\widetilde{V};\widetilde{F})$ is infinite-dimensional. This
enables us to get the following corollary.

\vskip 4mm

\no{\bf Corollary $2.2$} \ {\it Let
$(\widetilde{V};\widetilde{F})$ be a multi-vector space with
$\widetilde{V}=\bigcup\limits_{i=1}^k V_i,
\widetilde{F}=\bigcup\limits_{i=1}^k F_i$. Then
$(\widetilde{V};\widetilde{F})$ is finite-dimensional if and only
if for any integer $i, 1\leq i\leq k$, $(V_i; +_i,\cdot_i)$ is
finite-dimensional.}

\vskip 4mm

\no{\bf Theorem $2.3$} \ {\it For a finite-dimensional
multi-vector space $(\widetilde{V};\widetilde{F})$, any two bases
have the same number of vectors.}

\vskip 3mm

{\it Proof} \ Let $\widetilde{V}=\bigcup\limits_{i=1}^k V_i$ and
$\widetilde{F}=\bigcup\limits_{i=1}^k F_i$. The proof is by the
induction on $k$. For $k=1$, the assertion is true by Theorem $4$
of Chapter $2$ in $[3]$.

For the case of $k=2$, notice that by a result in linearly vector
space theory (see also [$3$]), for two subspaces $W_1,W_2$ of a
finite-dimensional vector space, if the basis of $W_1\bigcap W_2$
is $\{{\bf a}_1,{\bf a}_2,\cdots , {\bf a}_t\}$, then the basis of
$W_1\bigcup W_2$ is

$$\{{\bf a}_1,{\bf a}_2,\cdots , {\bf a}_t, {\bf b}_{t+1},{\bf b}_{t+2},\cdots ,{\bf b}_{dimW_1},
{\bf c}_{t+1},{\bf c}_{t+2}, \cdots ,{\bf c}_{dimW_2}\},$$

\no where, $\{{\bf a}_1,{\bf a}_2,\cdots , {\bf a}_t, {\bf
b}_{t+1},{\bf b}_{t+2},\cdots ,{\bf b}_{dimW_1}\}$ is a basis of
$W_1$ and $\{{\bf a}_1,{\bf a}_2,\cdots , {\bf a}_t,$ ${\bf
c}_{t+1},{\bf c}_{t+2},\cdots ,{\bf c}_{dimW_2}\}$ a basis of
$W_2$.

Whence, if $\widetilde{V} = W_1\bigcup W_2$ and
$\widetilde{F}=F_1\bigcup F_2$, then the basis of $\widetilde{V}$
is also

$$\{{\bf a}_1,{\bf a}_2,\cdots , {\bf a}_t, {\bf b}_{t+1},{\bf b}_{t+2},\cdots ,{\bf b}_{dimW_1},
{\bf c}_{t+1},{\bf c}_{t+2}, \cdots ,{\bf c}_{dimW_2}\}.$$

Assume the assertion is true for $k=l, l\geq 2$. Now we consider
the case of $k=l+1$. In this case, since

$$\widetilde{V}=(\bigcup\limits_{i=1}^l V_i)\bigcup V_{l+1}, \ \widetilde{F}
=(\bigcup\limits_{i=1}^l F_i)\bigcup F_{l+1},$$

\no by the induction assumption, we know that any two bases of the
multi-vector space $(\bigcup\limits_{i=1}^l
V_i;\bigcup\limits_{i=1}^l F_i)$ have the same number $p$ of
vectors. If the basis of $(\bigcup\limits_{i=1}^l V_i)\bigcap
V_{l+1}$ is $\{{\bf e}_1,{\bf e}_2,\cdots ,{\bf e}_n\}$, then the
basis of $\widetilde{V}$ is

$$\{{\bf e}_1,{\bf e}_2,\cdots ,{\bf e}_n, {\bf f}_{n+1},{\bf f}_{n+2},\cdots ,{\bf f}_{p},
{\bf g}_{n+1},{\bf g}_{n+2},\cdots ,{\bf g}_{dimV_{l+1}}\},$$

\no where $\{{\bf e}_1,{\bf e}_2,\cdots ,{\bf e}_n, {\bf
f}_{n+1},{\bf f}_{n+2},\cdots ,{\bf f}_{p}\}$ is a basis of
$(\bigcup\limits_{i=1}^l V_i; \bigcup\limits_{i=1}^l F_i)$ and
$\{{\bf e}_1,{\bf e}_2,\cdots ,{\bf e}_n,$ ${\bf g}_{n+1},{\bf
g}_{n+2},\cdots ,{\bf g}_{dimV_{l+1}}\}$ a basis of $V_{l+1}$.
Whence, the number of vectors in a basis of $\widetilde{V}$ is
$p+dimV_{l+1}-n$ for the case $n=l+1$.

Therefore, by the induction principle, we know the assertion is
true for any integer $k$.\quad\quad $\natural$

The number of a finite-dimensional multi-vector space
$\widetilde{V}$ is called its {\it dimension}, denoted by
$dim\widetilde{V}$.

\vskip 4mm

\no{\bf Theorem $2.4$}({\it dimensional formula}) \ {\it For a
multi-vector space $(\widetilde{V}; \widetilde{F})$ with
$\widetilde{V}=\bigcup\limits_{i=1}^k V_i$ and
$\widetilde{F}=\bigcup\limits_{i=1}^k F_i$, the dimension
$dim\widetilde{V}$ of $\widetilde{V}$ is}

$$dim\widetilde{V}=\sum\limits_{i=1}^k(-1)^{i-1}
\sum\limits_{\{i1,i2,\cdots ,ii\}\subset\{1,2,\cdots
,k\}}dim(V_{i1} \bigcap V_{i2}\bigcap\cdots\bigcap V_{ii}).$$

\vskip 3mm

{\it Proof} \ The proof is by induction on $k$. For $k=1$, the
formula is the trivial case of $dim\widetilde{V}=dimV_1$. for
$k=2$, the formula is

$$dim\widetilde{V}=dimV_1+dimV_2-dim(V_1\bigcap dimV_2),$$

\no which is true by Theorem $6$ of Chapter $2$ in $[3]$.

Now assume the formula is true for $k=n$. Consider the case of
$k=n+1$. According to the proof of Theorem $2.15$, we know that

\begin{eqnarray*}
dim\widetilde{V} &=& dim(\bigcup\limits_{i=1}^nV_i)+dimV_{n+1}
-dim((\bigcup\limits_{i=1}^nV_i)\bigcap V_{n+1})\\
 &=& dim(\bigcup\limits_{i=1}^nV_i)+dimV_{n+1}
-dim(\bigcup\limits_{i=1}^n(V_i\bigcap V_{n+1}))\\
&=& dimV_{n+1}+
\sum\limits_{i=1}^n(-1)^{i-1}\sum\limits_{\{i1,i2,\cdots,ii\}\subset\{1,2,\cdots
,n\}}dim(V_{i1}\bigcap V_{i2}\bigcap\cdots\bigcap V_{ii})\\
&+&
\sum\limits_{i=1}^n(-1)^{i-1}\sum\limits_{\{i1,i2,\cdots,ii\}\subset\{1,2,\cdots
,n\}}dim(V_{i1}\bigcap V_{i2}\bigcap\cdots\bigcap V_{ii}\bigcap
V_{n+1})\\
&=& \sum\limits_{i=1}^n(-1)^{i-1} \sum\limits_{\{i1,i2,\cdots
,ii\}\subset\{1,2,\cdots ,k\}}dim(V_{i1} \bigcap
V_{i2}\bigcap\cdots\bigcap V_{ii}).
\end{eqnarray*}

By the induction principle, we know this formula is true for any
integer $k$. \quad\quad $\natural$ \vskip 2mm

From Theorem $2.4$, we get the following additive formula for any
two multi-vector spaces.

\vskip 4mm

\no{\bf Corollary $2.3$}({\it additive formula}) \ {\it For any
two multi-vector spaces $\widetilde{V}_1,\widetilde{V}_2$,}

$$dim(\widetilde{V}_1\bigcup\widetilde{V}_2)=dim\widetilde{V}_1+
dim\widetilde{V}_2-dim(\widetilde{V}_1\bigcap\widetilde{V}_2).$$

\vskip 8mm

{\bf $3.$ Open problems for a multi-ring space}

\vskip 5mm

Notice that Theorem $2.3$ has told us there is a similar linear
theory for multi-vector spaces, but the situation is more complex.
Here, we present some open problems for further research.

\vskip 3mm

\no{\bf Problem $3.1$} \ {\it Similar to linear spaces, define
linear transformations on multi-vector spaces. Can we establish a
new matrix theory for linear transformations?}

\vskip 3mm

\no{\bf Problem $3.2$} \ {\it Whether a multi-vector space must be
a linear space?}

\vskip 2mm

\no{\bf Conjecture $A$} \ {\it There are non-linear multi-vector
spaces in multi-vector spaces.}\vskip 2mm

Based on Conjecture $A$, there is a fundamental problem for
multi-vector spaces.\vskip 2mm

\no{\bf Problem $3.3$} \ {\it Can we apply multi-vector spaces to
non-linear spaces?}

\vskip 8mm

{\bf References}

\vskip 5mm

\re{[1]}G.Birkhoff and S.Mac Lane, {\it A Survey of Modern
Algebra}, Macmillan Publishing Co., Inc, 1977.

\re{[2]}Daniel Deleanu, {\it A Dictionary of Smarandache
Mathematics}, Buxton University Press, London \& New York,2004.

\re{[3]}K.Hoffman and R.Kunze, {\it Linear Algebra} (Second
Edition), Prentice-Hall, Inc., Englewood Cliffs, New Jersey, 1971.

\re{[4]}L.F.Mao, On Algebraic Multi-Group Spaces, {\it eprint
arXiv: math/0510427}, 10/2005.

\re{[5]}L.F.Mao, {\it Automorphism Groups of Maps, Surfaces and
Smarandache Geometries}, American Research Press, 2005.

\re{[6]} F.Smarandache, Mixed noneuclidean geometries, {\it eprint
arXiv: math/0010119}, 10/2000.

\re{[7]}F.Smarandache, {\it A Unifying Field in Logics.
Neutrosopy: Neturosophic Probability, Set, and Logic}, American
research Press, Rehoboth, 1999.

\re{[8]}F.Smarandache, Neutrosophy, a new Branch of Philosophy,
{\it Multi-Valued Logic}, Vol.8, No.3(2002)(special issue on
Neutrosophy and Neutrosophic Logic), 297-384.

\re{[9]}F.Smarandache, A Unifying Field in Logic: Neutrosophic
Field, {\it Multi-Valued Logic}, Vol.8, No.3(2002)(special issue
on Neutrosophy and Neutrosophic Logic), 385-438.

\re{[10]]}W.B.Vasantha Kandasamy, {\it Bialgebraic structures and
Smarandache bialgebraic structures}, American Research Press,
2003.

\re{[11]}W.B.Vasantha Kandasamy and F.Smarandache, {\it Basic
Neutrosophic Algebraic Structures and Their Applications to Fuzzy
and Neutrosophic Models}, HEXIS, Church Rock, 2004.

\end{document}